\documentclass[11pt,twoside]{article}
\usepackage{amsfonts}
\begin{document}

\centerline{\bf The abc - conjecture is true }
 \centerline{\bf for at least $N(c), \hskip 4pt 1\leq N(c) < \varphi(c)/2$, partitions $a$, $b$ of $c$}
\par \bigskip \bigskip \centerline{Constantin M. Petridi,} \centerline{11 Apollonos
Street, 151 24 Maroussi (Athens), Greece, cpetridi@hotmail.com}
\par
\bigskip
\begin{abstract}
We prove that for any positive integer c there are at least N(c),
$1\leq N(c) < \varphi(c)/2$ representations of c as a sum of two
positive integers a, b, with no common divisor, such that the N(c)
radicals R(abc) are all greater than kc, where k an absolute
constant.
\end{abstract}
\bigskip
\par

 \textbf{Preliminaries.} Let
$c=q_1^{a_1}q_2^{a_2} \cdots q_{\omega}^{a_{\omega}}$, $q_i$
diferent primes, $a_i\geq 1$, be a positive integer. Denote its
radical $q_1q_2\cdots q_{\omega}$ by $R(c)$, and similarly $R(n)$
for any integer $n$. Consider the positive solutions of the
Diophantine equation $x+y=c$, $(x,y)=1$, $x< y$. Their number is
$\varphi (c)/2$ ($\varphi$ : Euler function). Denoting them in
some order by $a_i,b_i$, $1\leq i\leq \varphi (c)/2$, and listing
them with their respective radicals, one has
$$a_1+b_1=c \qquad R(a_1b_1c)$$
$$a_2+b_2=c \qquad R(a_2b_2c) \eqno (1)$$
$$\cdots\qquad\ \qquad\cdots $$
$$a_{\varphi (c) \over 2}+b_{\varphi (c)\over 2}=c \qquad
R(a_{\varphi (c) \over 2}b_{\varphi (c) \over 2}c).$$
\par
 Form the product of above radicals
$$G_c=\prod_{1\leq i \leq {\varphi(c) \over 2}}R(a_ib_ic).$$
\par
$G_c^{2/ \varphi(c)}$ is the geometric mean of the radicals.
\par
Throughout the paper $p$ designates prime numbers.
\par
\bigskip
\textbf{Definition.} We introduce a function $E_c(x)$ defined for
any real $x$ by $$E_c(x):=\biggl [  {1\over x}{c\over 1} \bigg ]
-\sum_{1\leq i\leq \omega} \bigg [ {1\over x}{c\over q_i} \bigg
]+\sum_{1\leq i,j\leq \omega}\bigg [ {1\over x}{c\over q_iq_j}
\bigg ] -\cdots +(-1)^{\omega} \bigg [ {1\over x}{c\over q_1
\cdots q_{\omega}} \bigg ], $$ which for positive integers $x$
plays a key role in our investigations. If the numbers in the
integral part brackets are integers the function reduces to
$$ {1\over x}{c\over 1} -\sum_{1\leq i\leq \omega} {1\over x}{c\over q_i}+\sum_{1\leq i,j\leq \omega}{1\over x}{c\over q_iq_j}
 -\cdots +(-1)^{\omega}  {1\over x}{c\over q_1 \cdots
q_{\omega}}={1\over x}\hskip 2pt {\varphi (c)}. \eqno (2) $$
\par
\bigskip
Following theorem gives an explicit expression for $G_c$.
\par
\textbf{Theorem 1.}$$G_c=R(c)^{\varphi (c) \over 2} \prod_{2\leq
p <c \atop (p,c)=1 }p^{E_c(p)} .$$
\par
 \textbf{Proof.} Consider \textit{all}
positive solutions of the Diophantine equation $x+y=c$, $x\leq y$,
and their corresponding radicals $R(xyc)$
$$1+(c-1)=c \qquad R(1(c-1)c)$$
$$2+(c-2)=c \qquad R(2(c-2)c) \eqno (3)$$
$$\cdots \qquad\qquad \qquad\qquad \cdots$$
$$\bigg [ {c \over 2} \bigg ]+ \bigg ( c- \bigg [ {c \over 2} \bigg ] \bigg )=c
\qquad R \bigg ( \bigg [{c \over 2} \bigg ]\bigg [ c- {c \over 2}
\bigg ]c \bigg ).$$ Equalities (3) comprise equalities (1), but
include also those for which $(x,y)>1$. That all $q_i$ appear
$\varphi (c) \over 2$  times in $G_c$ is obvious. Consequently,
so does their product $R(c)$. For primes $p \neq q_i$, we apply
the inclusion-exclusion principle. As all numbers $<c$ do occur in
the equalities (3), the number of times $p$ appears in the
radicals (3) is $\big [ {c\over p} \big ]$. The number of times
$pq_i$ appears in the radicals (3) is $\sum_{i=1}^{\omega} \big [
{1\over {pq_i}} \big ]$. The number of times $pq_iq_j$ appears in
the radicals (3) is $\sum_{i,j=1}^{\omega} \big [ {1\over
{pq_iq_j}} \big ], $ e.t.c. Inserting these values in the
inclusion-exclusion formula we get for the total number of times
$p$ appears in the radicals (3)
$$\biggl [  {1\over p}{c\over 1} \bigg
] -\sum_{1\leq i\leq \omega} \bigg [ {1\over p}{c\over q_i} \bigg
]+\sum_{1\leq i,j\leq \omega}\bigg [ {1\over p}{c\over q_iq_j}
\bigg ] -\cdots +(-1)^{\omega} \bigg [ {1\over p}{c\over q_1
\cdots q_{\omega}} \bigg ]=E_c(p), $$ as stated in the theorem.
\par
\bigskip
\textbf{Corollary 1. } If $d$ denotes the divisors of $c$, then
$$\prod_{d|c}G_d=q_1^{\Theta
(q_1)} \ldots q_{\omega}^{\Theta (q_{\omega})} \prod_{2\leq p <c
\atop {(p,c)=1}}p^{[ {c\over p} ]}, $$ where
$$\Theta (q_i)= \sum_{d\equiv 0(q_i) \atop d|c}{\varphi (d) \over
2}+\sum_{d\neq 0(q_i) \atop d|c}E_d(q_i).$$
\textbf{Proof.}
Applying Theorem 1 to the divisors $d$ of $c$ one has
$$G_d=R(d)^{\varphi (d) \over 2} \prod_{2 \leq p < d \atop
(p,d)=1}p^{E_d(p)}.$$ Multiplying over all $\tau (c)$ divisors and
using the, easily established, fact that
$$\prod_{x+y=c \atop x\leq y}R(xyc)=R(c)^{[{c\over 2}
]}\prod_{2\leq p <c \atop (p,c)=1}p^{[{c\over p}]},$$ gives for
$\hskip 4pt \prod_{d|c}G_d \hskip 4pt$ the result.
\par
The corollary will not be used in the sequel.

\par
\bigskip
We shall now prove certain Lemmas regarding the function
$E_c(p)$. We shall also list, without proof, some well known facts
from the elementary theory of primes so as not to interrupt the
main body of the proof. Absolute constants will be denoted by
$k_i$, indexed in the order they first appear.
\par
\textbf{Lemma 1.} For $x>0$ (actually for any $x\neq 0$)
$${\varphi (c) \over x}-2^{\omega -1} < E_c(x) <
{\varphi (c) \over x}+2^{\omega -1}. \eqno (4) $$

\textbf{Proof.} By definition one has
$$ {1 \over x}{c \over 1}-{\omega \choose 0} < \bigg [{1 \over x}{c \over 1} \bigg ] \leq {1 \over x}{c \over
1}$$
$$ -\sum_{1\leq i \leq \omega} {1\over x}{c\over
q_i} \leq -\sum_{1\leq i\leq \omega} \bigg [ {1\over x}{c\over
q_i} \bigg ] < {\omega \choose 1}- \sum_{1\leq i\leq \omega}
{1\over x}{c\over q_i}$$
$$\sum_{1\leq i,j\leq \omega} {1\over x}{c\over
q_iq_j}-{\omega \choose 2} < \sum_{1\leq i,j\leq \omega}\bigg [
{1\over x}{c\over q_iq_j} \bigg ] \leq  \sum_{1\leq i,j\leq
\omega} {1\over x}{c\over q_iq_j} $$
$$\qquad \cdots \qquad \qquad \qquad \cdots \qquad \qquad \qquad \cdots $$
$$  {1\over x}{c\over
q_1 \cdots q_{\omega}} -{\omega \choose \omega}  < (-1)^{\omega}
\bigg [ {1\over x}{c\over q_1 \cdots q_{\omega}} \bigg ] \leq
{1\over x}{c\over q_1 \cdots q_{\omega}}, \qquad \omega \equiv
0(2)$$
$$  -{1\over x}{c\over
q_1 \cdots q_{\omega}} \leq (-1)^{\omega} \bigg [ {1\over
x}{c\over q_1 \cdots q_{\omega}} \bigg ] < {\omega \choose
\omega}- {1\over x}{c\over q_1 \cdots q_{\omega}}, \qquad \omega
\equiv 1(2).$$ Adding term-wise above inequalities, we have
$${c \over x} \bigg \{ 1- \sum{1 \over q_i} + \sum{1 \over
q_iq_j}- \cdots \bigg \}-\sum_{\nu\equiv 0(2)}{\omega \choose
{\nu}} < E_c(x) < {c \over x} \bigg \{ 1- \sum{1 \over q_i} +
\sum{1 \over q_iq_j}- \cdots \bigg \}+\sum_{\nu \equiv
1(2)}{\omega \choose {\nu}}.$$ Considering that
$$\sum_{\nu=0(2)}{ \omega \choose {\nu}}=\sum_{\nu=1(2)}{ \omega \choose
{\nu}}=2^{\omega -1}, $$ we have by (2), as required,
$${\varphi (c) \over x}-2^{\omega -1} < E_c(x) < {\varphi (c) \over x}+2^{\omega
-1}. $$
\par
\bigskip
\textbf{Lemma 2.}
$$   E_c(x) \quad > \quad
\left\{\begin{tabular}{c@{\quad}l}
                       ${{\varphi (c)} \over {x}}-2^{\omega -1}$
                       & \quad for \quad $0 < x < {{\varphi (c)} \over {2^{\omega-1}}}$ \\
                       $1$ & \quad for \quad ${{\varphi (c)} \over
                       {2^{\omega-1}}} < x < c$.

                                              \end{tabular} \right.        $$

\textbf{Proof.} The expression ${{\varphi (c)} \over
{x}}-2^{\omega -1}$ is positive for all $0<x<{{\varphi (c)} \over
{2^{\omega -1}}}$. For $x>{{\varphi (c)} \over {2^{\omega -1}}}$
the lowest limit of $E_c(x)$ is $1$, since $x<c$.
\par
\bigskip
\textbf{Lemma 3.} If one (or more) of the prime factors $q_i$, has
exponent $a_i \geq 2$ then $$E_c(q_i)= {{\varphi (c)} \over
q_i}.$$ \textbf{Proof.} Writing $q_i$ instead of $x$ in (2), and
renaming the running indexes, one has
$$E_c(q_i)=\bigg [  {1\over q_i}{c\over 1} \bigg ] -\sum_{1\leq
j\leq \omega} \bigg [ {1\over q_i}{c\over q_j} \bigg ]+\sum_{1\leq
j,k\leq \omega}\bigg [ {1\over q_i}{c\over q_jq_k} \bigg ] -\cdots
+(-1)^{\omega} \bigg [ {1\over q_i}{c\over q_1 \cdots q_{\omega}}
\bigg ]. $$ Since by supposition the exponent of $q_i$ is at least
$\geq 2$ the numbers within the integral part brackets are
integers so that we can skip the brackets. This gives by (2)
$$ {1\over q_i}{c\over 1} -\sum_{1\leq j\leq \omega} {1\over q_i}{c\over q_j}+\sum_{1\leq j,k\leq \omega}{1\over q_i}{c\over q_jq_k}
 -\cdots +(-1)^{\omega}  {1\over q_i}{c\over q_1 \cdots
q_{\omega}}={\varphi(c) \over q_i} , $$ as stated.
\par
\bigskip
\textbf{Lemma 4.} If one (or more) of the prime factors $q_i$ has
exponent $a_i=1$, then
$$E_c(q_i)={\varphi (c) \over {q_i -1}}-E_{c / q_i}(q_i).$$
\textbf{Proof.} For convenience, putting $\overline{c}=c/q_i$,
renaming indexes as above, and writing on the left hand side the
terms of the $E_c(q_i)$- function vertically, we have following
equalities
$$\bigg [ {1\over q_i}{c\over 1} \bigg ] = \bigg [ {1\over
q_i}{q_i \hskip 2pt \overline{c}\over 1} \bigg ] $$
$$- \sum_{1\leq i \leq \omega} \bigg [ {1\over q_i}{c\over q_j}
\bigg ] = - \bigg [ {1\over q_i}{q_i \hskip 2pt \overline{c} \over
q_i} \bigg ]- \sum_{2\leq j \leq \omega} \bigg [ {1\over q_i}{q_i
\hskip 2pt \overline{c}\over q_j} \bigg ]$$
$$+ \sum_{1\leq j,k \leq \omega} \bigg [ {1\over q_i}{c \over {q_jq_k}} \bigg
]= +\sum_{2\leq j \leq \omega} \bigg [ {1\over q_i}{q_i \hskip 2pt
\overline{c} \over {q_i q_j}} \bigg ]+\sum_{2\leq j,k \leq \omega}
\bigg [ {1\over q_i}{q_i \hskip 2pt \overline{c} \over {q_jq_k}}
\bigg ]$$
$$\qquad \ldots \qquad \qquad \qquad \ldots \qquad \qquad \qquad \ldots \qquad$$
$$(-1)^{\omega -1} \sum_{1\leq j_{\nu } \leq {\omega -1} }\bigg [ {1\over
q_i}{c\over {q_{j_1} \ldots q_{j_{\omega -1}}}} \bigg ] =
(-1)^{\omega -1} \sum_{2\leq j_{\nu } \leq {\omega -1} }\bigg [
{1\over q_i}{q_i \hskip 2pt \overline{c}\over {q_iq_{j_1} \ldots
q_{j_{\omega -1}}} }  \bigg ]+(-1)^{\omega -1} \bigg [ {1\over
q_i} {q_i \hskip 2pt \overline{c}\over {q_2 \ldots q_{\omega}} }
\bigg ]$$
$$(-1)^{\omega } \bigg [ {1\over q_i}{c \over {q_1q_2 \ldots
q_{\omega}}} \bigg ] = (-1)^{\omega } \bigg [ {1\over q_i}{q_i
\hskip 1pt \overline{c}\over {q_iq_2 \ldots q_{\omega}}} \bigg ]$$
\par
 Adding these equalities, the sum of the left hand side terms is, as
said, $E_c(q_i)$.
\par
The right hand side is equal to
$$ \bigg \{ - \bigg [ {1\over q_i}{q_i \hskip 2pt \overline{c} \over q_i}
\bigg ]+\sum_{2\leq j \leq \omega} \bigg [ {1\over q_i}{q_i \hskip
2pt \overline{c} \over {q_i q_j}} \bigg ]-\cdots+(-1)^{\omega -1}
\sum_{2\leq j_{\nu } \leq {\omega -1} }\bigg [ {1\over q_i}{q_i
\hskip 2pt \overline{c}\over {q_iq_{j_1} \ldots q_{j_{\omega -2}}}
} \bigg ]+(-1)^{\omega } \bigg [ {1\over q_i}{q_i \hskip 2pt
\overline{c}\over {q_iq_2 \ldots q_{\omega}}} \bigg ] \bigg \}$$
$$+ \bigg \{\bigg [ {1\over
q_i}{q_i \hskip 2pt \overline{c}\over 1} \bigg ]-\sum_{2\leq j\leq
\omega} \bigg [ {1\over q_i}{q_i \hskip 2pt \overline{c}\over q_j}
\bigg ]+\sum_{2\leq j,k \leq \omega} \bigg [ {1\over q_i}{q_i
\hskip 2pt \overline{c} \over {q_jq_k}} \bigg
]-\cdots+(-1)^{\omega -1} \bigg [ {1\over q_i} {q_i \hskip 2pt
\overline{c}\over {q_2 \ldots q_{\omega}} } \bigg ] \bigg \}. $$
Cancelling $q_i$ within the brackets, the first brace is clearly
equal to $-E_{\overline{c}}(q_i)$. In the second brace, as all
$q_j,q_k$, $2\leq j,k \leq \omega $ divide $\overline{c}$, the
numbers within the integral part brackets are integers. We
therefore can write it
$$\overline{c} \hskip 2pt  \bigg ( 1- \sum_{2\leq j \leq \omega}{1\over q_j}
+ {\sum_{2\leq j,k \leq \omega}}{1\over
{q_jq_k}}-\cdots+(-1)^{\omega -1}{1\over{q_2 \ldots q_{\omega}}}
\bigg )= {\overline{c}} \hskip 2pt {\prod_{2\leq j\leq \omega}}
\bigg ( 1-{1\over q_j} \bigg )=  \varphi (\overline{c}) .$$
Substituting the found values in the sum of above equalities and
noting that $\varphi(\overline{c})={\varphi(c)}/{(q-1)}$, since
$(c, \overline{c})=1$, we finally get
$$E_c(q_i)={\varphi(c)\over {q_i-1}}-E_{c/q_i}(q_i),$$
which proves the Lemma.
\par
\bigskip
\textbf{Lemma 5.} $$e^{k_1c} < \prod_{2\leq p\leq c}p < e^{k_2c}.
$$ \textbf{Proof.} This is the multiplicative form of
Tchebycheff's estimate for $\sum_{2\leq p\leq c}log \hskip 2pt
p$, namely,
$$k_1c<\sum_{2\leq p\leq c}log \hskip 2pt p < k_2c ,$$ for $c\geq 2$,
with $k_1,k_2$ positive absolute constants.
\par
\bigskip
\textbf{Lemma 6.}$$e^{-k_3}c < \prod_{2\leq p\leq c}p^{1\over p} <
e^{k_3}c.$$
 \textbf{Proof.} This is the multiplicative form of
Tchebycheff's estimate for $\sum_{2\leq p\leq c}{1\over p}\hskip
2pt log \hskip 2pt p$, namely,
$$logc-k_3<\sum_{2\leq p\leq c}{1\over p}\hskip 2pt log \hskip 2pt p < logc+k_3 ,$$ for $c\geq 2$, with $k_3$
a positive absolute constant.
\par
\bigskip
We now state a result which gives a lower bound for the geometric
mean $G_c^{2 \over \varphi (c)}$ in terms of the prime factors
$q_i$ of $c$ and their respective exponents $a_i$.
\par
\textbf{Theorem 2.}$$G_c^{2\over \varphi (c)} > k_4 \hskip 2pt
\prod_{1\leq i\leq \omega} \hskip 2pt q_i^{2a_i-1-{2 \over \varphi
(c)}E_c(q_i)}\bigg ( {{q_i-1}\over 2} \bigg )^2  ,$$ where $k_4$
 a positive absolute constant.
\par
 \textbf{Proof.} We transform the expression given for $G_c$ in
 Theorem 1, as follows:
 $$G_c=\prod_{1\leq i \leq \omega}q_i^{\varphi (c) \over 2} \hskip 2pt
 \prod_{2\leq p<c \atop (p,c)=1} p^{E_c(p)} \qquad \qquad \qquad \qquad \qquad \qquad \qquad \qquad$$
 $$\qquad ={{\prod_{1\leq i \leq \omega}q_i^{\varphi (c) \over 2}} {\prod_{1\leq i \leq \omega}
 q_i^{-E_c(q_i)}}} \hskip 2pt {\prod_{2\leq p \leq c}p^{E_c(p)}}. \hskip 2pt
 \qquad \qquad \qquad \qquad \qquad \quad$$
Joining the first two products into one and splitting the third
product as indicated we have
$$G_c ={\prod_{1\leq i \leq \omega}q_i^{{\varphi
(c) \over 2}-E_c(q_i)}} \hskip 2pt{\prod_{2\leq p< {\varphi (c)
\over 2^{\omega -1}}}p^{E_c(p)} } \hskip 2pt \prod_{{\varphi (c)
\over 2^{\omega -1}}<p<c}p^{E_c(p)}. \quad \qquad \qquad$$
Applying Lemma 2 to the second and third product and splitting
the products in an obvious way, we get successively
$$ G_c > {\prod_{1\leq i \leq \omega}q_i^{{\varphi
(c) \over 2}-E_c(q_i)}} \hskip 2pt {\prod_{2\leq p<{\varphi (c)
\over 2^{\omega -1}}}p^{{\varphi (c) \over p}-2^{\omega -1}}}
\hskip 2pt \prod_{{\varphi (c) \over 2^{\omega -1}}<p<c}p \qquad
\qquad \qquad \qquad \qquad \qquad$$
$$\qquad >{\prod_{1\leq i \leq \omega}q_i^{{\varphi
(c) \over 2}-E_c(q_i)}} \hskip 2pt \bigg \{\prod_{2\leq p<{\varphi
(c) \over 2^{\omega -1}}}p^{1\over p} \bigg \}^{\varphi (c)} \bigg
\{\prod_{2\leq p<{\varphi (c) \over 2^{\omega -1}}}p \bigg
\}^{-2^{\omega -1}} \hskip 2pt  \bigg \{\prod_{2\leq p<{\varphi
(c) \over 2^{\omega -1}}}p \bigg \}^{-1} \hskip 2pt {\prod_{2\leq
p<c}p}. \hskip 8pt$$ Joining the third and the fourth product into
one, we have
$$G_c >{\prod_{1\leq i \leq \omega}q_i^{{\varphi
(c) \over 2}-E_c(q_i)}} \hskip 2pt \bigg \{\prod_{2\leq p<{\varphi
(c) \over 2^{\omega -1}}}p^{1\over p} \bigg \}^{\varphi (c)}
 \hskip 2pt \bigg \{\prod_{2\leq p<{\varphi
(c) \over 2^{\omega -1}}}p \bigg \}^{-(2^{\omega -1}+1)} \hskip
2pt {\prod_{2\leq p<c}p}. \qquad \quad \quad$$ Applying Lemma 6 to
the second product, Lemma 5 to the third and fourth product, we
have
$$G_c >{\prod_{1\leq i \leq \omega}q_i^{{\varphi
(c) \over 2}-E_c(q_i)}} \hskip 2pt \bigg \{e^{-k_3}{\varphi (c)
\over {2^{\omega -1}}} \hskip 2pt \bigg \}^{\varphi (c)} \hskip
2pt \bigg \{e^{k_2 {\varphi (c) \over {2^{\omega -1}}}} \bigg
\}^{-(2^{\omega -1}+1)}e^{k_1c}. \qquad \qquad$$ Summing the
exponents of $e$, we have
$$G_c > \bigg ( {\prod_{1\leq i \leq \omega}q_i^{{\varphi
(c) \over 2}-E_c(q_i)}} \bigg ) \hskip 2pt \bigg ( {\varphi (c)
\over 2^{\omega -1}} \bigg )^{\varphi (c)}e^{-k_3 \varphi (c) +
k_1 c - k_2 \varphi (c) -k_2 {\varphi (c) \over 2^{\omega -1}}}.
\qquad \qquad \qquad \quad$$ Raising this inequality to the power
$2 \over \varphi (c)$, we get for the geometric mean $G_c^{2
\over \varphi (c)}$
$$G_c^{2 \over \varphi (c)} > \bigg ( {\prod_{1\leq i \leq \omega}
q_i^{1- {2\over \varphi (c)}E_c(q_i)}} \bigg ) {\bigg ( {{\varphi
(c)} \over 2^{\omega -1}}} \bigg )^2 e^{-2k_3+2k_1{c\over \varphi
(c)}-2k_2 (1+{1\over {2^{\omega -1}}})}.  \quad \quad \qquad
\qquad \eqno (5) $$ Evaluating the second parenthesis we have
 $${\bigg ( {{\varphi
(c)} \over 2^{\omega -1}}} \bigg )^2 = 4{\prod_{1\leq i \leq
\omega}q_i^{2a_i-2}} \hskip 2pt \bigg ( {{q_i-1} \over 2} \bigg
)^2. $$ On the other hand, since ${c \over \varphi (c)}>1 $ and
$1+{1\over 2^{\omega -1}} \leq 2$, the entire exponent of $e$
appearing in (4) is $> -2k_3+2k_1-4k_2$. Setting
$k_4=4e^{2k_1-4k_2-2k_3}$ as a new absolute constant and
substituting in (5) we finally get
$$G_c^{2\over \varphi (c)} > k_4 \hskip 2pt
\prod_{1\leq i\leq \omega} \hskip 2pt q_i^{2a_i-1-{2 \over \varphi
(c)}E_c(q_i)} \bigg ( {{q_i-1}\over 2} \bigg )^2 ,$$ which
completes the proof.
\par
\bigskip
\textbf{Theorem 3.}
$$G_c^{2\over {\varphi (c)}}> k \hskip 2pt c \hskip 2pt ,$$ where $k$ a positive
absolute constant.
\par
 \textbf{Proof.} Denote the expressions within the product of
 Theorem 2 by $F(q_i,a_i)$, $1\leq i\leq \omega$.
\par
For $\underline{a_i=1}$ by application of Lemma 4 and taking into
account that for $x>0$, $E_c(x)\geq 0$ we have
$$F(q_i,a_i)=q_i^{2a_i-1-{2\over \varphi (c)}E_c(q_i)} \bigg (
{{q_i-1}\over 2 }\bigg )^2$$
$$\qquad\quad = q_i^{1-{2\over \varphi (c)}E_c(q_i)}\bigg (
{{q_i-1}\over 2 }\bigg )^2 \hskip 8pt$$
$$\qquad\qquad\qquad = q_i^{{{q_i-3}\over {q_i-1}}+{2 \over \varphi
(c)}E_{c/q_i}(q_i)}\bigg ( {{q_i-1}\over 2 }\bigg )^2 \hskip
16pt$$
$$\quad\quad  > q_i^{{q_i-3}\over {q_i-1}}\bigg ( {{q_i-1}\over 2 }\bigg )^2. \hskip 25pt$$
\par
Elementary considerations give for $q_i \geq 5$
$$F(q_i,1)>q_i^{{q_i-3}\over {q_i-1}}\bigg ( {{q_i-1}\over 2 }\bigg
)^2> \bigg ( {1\over 1} \bigg )q_i.$$ The exceptions are $q_i=2$,
resp. $3$ for which
$$F(2,1) > \bigg ( {1\over 16} \hskip 2pt \bigg ) 2$$
$$F(3,1) > \bigg ( {1\over 3} \hskip 2pt  \bigg ) 3.$$
For $\underline{a_i\geq 2}$ we have by application of Lemma 3
$$F(q_i,a_i)= q_i^{2a_i-1-{2\over \varphi (c)}E_c(q_i)}\bigg (
{{q_i-1}\over 2 }\bigg )^2$$
$$\qquad\quad =q_i^{2a_i-1-{2\over q_i}}\bigg (
{{q_i-1}\over 2 }\bigg )^2. \hskip 18pt$$
\par
Elementary considerations give for $q_i\geq 3$
$$F(q_i,a_i) > q_i^{a_i}.$$
\par
The exception is $q_i=2$ for which
$$F(2,a_i) > \bigg ( {1\over 4} \bigg ) 2^{a_i}.$$
Summarizing, above inequalities take one of following six forms,
depending on the, mutually exclusive, combinations of the
exceptional factors $F(2,1),F(3,1),F(2,a_i)$ that may occur in
the product $\prod_{1\leq i\leq \omega} F(q_i,a_i)$
$$\begin{tabular}{lllll}
$G_c^{2\over \varphi (c)}$& $>$& $k_4 \hskip 2pt c$ & exceptional
factor &\quad ---\qquad \\
$G_c^{2\over \varphi (c)}$&$ >$&$ k_4 \hskip 2pt {1\over 16}\hskip
2pt c$ & exceptional factor &$ F(2,1)$\\
$G_c^{2\over \varphi (c)}$&$>$&$ k_4 \hskip 2pt {1\over 3}\hskip
2pt c $& exceptional factor&$F(3,1)$\\
$G_c^{2\over \varphi (c)}$&$ >$&$ k_4 \hskip 2pt {1\over 16}\hskip
2pt {1\over 3} \hskip 2pt c$& exceptional factor&$
F(2,1),F(3,1)$\\
$G_c^{2\over \varphi (c)}$&$ >$&$ k_4 \hskip 2pt {1\over 4}\hskip
2pt c $& exceptional factor&$ F(2,a_i)$\\
$G_c^{2\over \varphi (c)}$&$ >$&$ k_4 \hskip 2pt {1\over 4}\hskip
2pt {1\over 3} \hskip 2pt c $& exceptional factor&$
F(2,a_i),F(3,1)$
\end{tabular}$$
 It results that all cases are included
in $$G_c^{2\over \varphi (c)} > k \hskip 2pt c,$$ where $k$ an
absolute constant $(= {1\over 48}k_4)$, as stated in the Theorem.
\par
\bigskip
\textbf{Theorem 4.} For any positive integer $c$ there is an
integer $N=N(c)$, $1\leq N < {\varphi (c) \over 2}$ such that at
least $N$ of the radicals figuring in (1), denoted generally by
$R(abc)$ satisfy
$$R(abc)>k \hskip 2pt c,$$ where $k$ the absolute constant of
Theorem 3.
\par
\textbf{Proof.} Arrange the ${\varphi(c)} \over 2$ radicals in (1)
in ascending order
$$R_{i_1},R_{i_2}, \ldots R_{i_{\varphi (c) \over 2}} .$$
 $G_c^{2\over \varphi(c)}$ is their geometric mean and hence
$$R_{i_1} < G_c^{2\over \varphi(c)} < R_{i_{\varphi (c) \over
2}}.$$ Let $N=N(c)$ be the number of radicals $> G_c^{2\over
\varphi(c)}$, $1\leq N < {\varphi(c) \over 2}$. Since by Theorem
3, $\hskip 2pt G_c^{2\over \varphi(c)}
> k \hskip 2pt c \hskip 4pt$ it follows a fortiori that for these $N$
radicals we have $$R(abc) > k \hskip 2pt c,$$ as was to be proved.
\par
\bigskip
\textbf{Comment.} For the radicals $>G_c^{2\over \varphi(c)},
\hskip 2pt$ the assertion of Theorem 4 is stronger than the
$abc$-conjecture
$$R(abc)>k(\varepsilon)^{-{1\over {1+\varepsilon}}} \hskip 2pt
c^{1\over {1+\varepsilon}}, \quad \textrm{for any} \hskip 6pt
\varepsilon
> 0,
$$ since the exponent of $c$ is $1$ and $k$ is absolute.
\par
It is weaker, however, than the $abc$-conjecture because it does
not say anything regarding the radicals $<G_c^{2\over
\varphi(c)}$.
\par
\bigskip
\textbf{Acknowledgment.} I am indebted to Peter Krikelis of the
University of Athens for his unfailing assistance.

\end{document}